\documentclass[10pt]{amsart}

\usepackage[latin1]{inputenc}
\usepackage[english]{babel}
\usepackage{indentfirst}
\usepackage{amssymb}
\usepackage{amsthm}
\usepackage{amsmath}

\newcommand{\sub}{\subseteq}
\def\cC{{{\mathcal C}}}
\def\cD{{{\mathcal D}}}
\def\cK{{{\mathcal K}}}
\def\N{\mathbb N}
\newcommand{\Ba}{{\rm Ba}}
\newcommand{\Bo}{{\rm Bo}}

\newcommand{\erre}{\mathbb{R}}
\def\epsilon{\varepsilon}

\newtheorem{theo}{theorem}[section]

\newtheorem{problem}[theo]{Problem}

\numberwithin{equation}{section}

\title{Open problems in Banach spaces and measure theory}

\author{Jos\'{e} Rodr\'{i}guez}
\address{Dpto. de Ingenier\'{i}a y Tecnolog\'{i}a de Computadores,
Facultad de Inform\'{a}tica, Universidad de Murcia, 30100 Espinardo (Murcia), Spain} \email{joserr@um.es}



\thanks{Research partially supported by {\em Ministerio de Econom\'{i}a y Competitividad} and {\em FEDER} (project MTM2014-54182-P).
This work was also partially supported by the research project 19275/PI/14 funded by {\em Fundaci\'{o}n S\'{e}neca - Agencia de Ciencia y Tecnolog\'{i}a 
de la Regi\'{o}n de Murcia} within the framework of {\em PCTIRM 2011-2014}.}

\begin{document}

\begin{abstract}
We collect several open questions in Banach spaces, mostly related to measure theoretic aspects of the theory. 
The problems are divided into five categories: miscellaneous problems in Banach spaces
(non-separable $L^p$ spaces, compactness in Banach spaces, $w^*$-null sequences in dual spaces),
measurability in Banach spaces (Baire and Borel $\sigma$-algebras, measurable selectors), 
vector integration (Riemann, Pettis and McShane integrals), vector measures (range and associated $L^1$ spaces) and
Lebesgue-Bochner spaces (topological and structural properties, scalar convergence).
\end{abstract}

\maketitle

\section*{Introduction}

The interaction between Banach space theory and measure theory has provided truly important results
in functional analysis. The topic reached its maturity in the seventies and eighties thanks to the contributions of outstanding mathematicians like 
J.~Bourgain, G.A.~Edgar, D.H.~Fremlin, I.~Namioka, R.R.~Phelps, C.~Stegall, M.~Talagrand, etc.
For detailed information we refer to the monographs \cite{bou-J,cem-men,die-uhl-J,tal,van}. 
Among the most celebrated results in this context, one finds the characterization of Asplund spaces
as those whose dual has the Radon-Nikod\'{y}m property, and the characterization of Banach spaces not
containing~$\ell^1$ as those whose dual has the weak Radon-Nikod\'{y}m property. 
Regarding applications to other areas, we can say, for instance, that the techniques of
set-valued integration in Banach spaces (initiated by R.J.~Aumann and G.~Debreu) 
have been used successfully in mathematical economics, control theory, game theory, etc. (see e.g. \cite{aub-fra,kle-tho}).

In this expository paper we collect some open questions (mostly) related to measure theoretic aspects of Banach space theory. 
The choice of problems is conditioned by our own research interests (and shortcomings) and we do not pretend 
to be exhaustive. The questions are organized in several sections:

\begin{itemize}
\item[1.] Miscellaneous problems in {\em Banach spaces}: non-separable $L^p$ spaces, compactness in Banach spaces, $w^*$-null sequences in dual spaces.
\item[2.] Problems in {\em measurability in Banach spaces}: Baire and Borel $\sigma$-algebras, measurable selectors.
\item[3.] Problems in {\em vector integration}: Riemann, Pettis and McShane integrals.
\item[4.] Problems in {\em vector measures}: range and associated $L^1$ spaces.
\item[5.] Problems in {\em Lebesgue-Bochner spaces}: topological and structural properties, scalar convergence.
\end{itemize}

We hope that this paper is a useful source of information as well as inspiration for researchers in Banach space
theory and measure theory.

\subsection*{Notation and terminology}

Our standard references are \cite{alb-kal,fab-ultimo} (Banach spaces) and \cite{die-uhl-J} (vector measures).
All our topological spaces are Hausdorff and we only consider real Banach spaces. 
An {\em operator} between Banach spaces is a linear continuous map. 
Given a Banach space~$X$, its norm is denoted by $\|\cdot\|_X$ or simply
$\|\cdot\|$. We write $B_X=\{x\in X:\|x\|\leq 1\}$ (the closed unit ball of~$X$) and $X^*$
denotes the topological dual of~$X$. The weak topology on~$X$ and the weak$^*$ topology on~$X^*$ 
are denoted by $w$ and~$w^*$, respectively. 
A {\em subspace} of~$X$ is a closed linear subspace. The subspace generated by a set~$A\sub X$ is denoted
by $\overline{{\rm span}}(A)$. Given another Banach space~$Y$, we say that $X$ contains~$Y$ or that $X$ contains a copy of~$Y$
(and we write $X \supseteq Y$) if there is a subspace of~$X$ which is isomorphic to~$Y$.
Recall that $X$ is said to be {\em Asplund} if every separable subspace of~$X$ has separable dual
(equivalently, $X^*$ has the Radon-Nikod\'{y}m property).
A set $\Gamma \sub X^*$ is said to be {\em total} if it separates the points of~$X$, i.e. for every $x\in X\setminus \{0\}$ 
there is $x^*\in \Gamma$ such that $x^*(x)\neq 0$. The topology on~$X$ of pointwise convergence on a total set~$\Gamma \sub X^*$ 
is denoted by $\sigma(X,\Gamma)$. 

The unit interval $[0,1]$ is equipped with the Lebesgue measure~$\lambda$ on the $\sigma$-algebra
of all Lebesgue measurable subsets. 
The symbols~$\mathfrak{c}$ and~$\omega_1$ denote the cardinality of the continuum and the first
uncountable ordinal, respectively. 
The {\em density character} of a topological space~$T$ is the cardinal
${\rm dens}(T)=\min\{|A|: A \sub T, \, \overline{A}=T\}$, where $|A|$ denotes the cardinality of the set~$A$.
We use the term ``compact'' as an abbreviation of ``compact topological space''.

\section{Miscellaneous problems in Banach spaces}\label{section:Banach}

\subsection{Non-separable $L^p$ spaces}\label{subsection:Lp-nosetoble}

Several aspects of the structure of non-separable $L^p$ spaces
are still unclear. In the separable case, it is known that $L^1[0,1]$ has no unconditional basis and, 
moreover, it is not a subspace of a Banach space with unconditional basis. 
On the other hand, $L^p[0,1]$ has unconditional basis for every $1<p<\infty$. 
The concept of unconditional basis admits a natural extension to the non-separable case (see e.g. \cite[Section~7.3]{fab-alt-JJ}). 
Enflo and Rosenthal~\cite{enf-ros} proved that $L^p(\mu)$ (for $1<p<\infty$, $p\neq 2$, and a finite measure~$\mu$) is not a subspace
of a Banach space with unconditional basis if ${\rm dens}(L^p(\mu))\geq \aleph_\omega$. They also conjectured that
the same conclusion holds whenever $L^p(\mu)$ is not separable.
The assumption on the density character is needed to use a combinatorial lemma that works only for ``big'' cardinals.
It is worth mentioning that a weakening of that lemma, for arbitrary uncountable sets, has been recently used 
in~\cite{avi-mar} to study certain extension operators between spaces of continuous functions.

\begin{problem}[Enflo-Rosenthal, \cite{enf-ros}]\label{problem:EnfloRosenthal}
Let $1<p<\infty$, $p\neq 2$, and let $\mu$ be a finite measure such that $L^p(\mu)$ is not separable. Can $L^p(\mu)$ be a subspace
of a Banach space having unconditional basis?
\end{problem}

In an attempt to attack the previous problem, Johnson and Schechtman~\cite{joh-sch} got interesting results
on non-separable $L^p$ spaces. A simple application of Pitt's theorem shows that, for $2<p<\infty$, the space
$\ell^p(\omega_1)$ is not a subspace of~$L^p(\mu)$ for any finite measure~$\mu$. It is more complicated
to prove the same statement for the range $1<p<2$,~\cite{enf-ros}. On the other hand, 
if $1<p<2$ and $X$ is a subspace of an~$L^p$ space (over an arbitrary non-negative measure)
such that $X\not\supseteq \ell^p(\omega_1)$, then $X \sub L^p(\mu)$ for some finite measure~$\mu$,~\cite{joh-sch}.
It is conjectured that this result also holds true for $2<p<\infty$. 

\begin{problem}[Johnson-Schechtman, \cite{joh-sch}]\label{problem:JS}
Let $2<p<\infty$ and let $X$ be a subspace of an~$L^p$ space such that $X\not\supseteq \ell^p(\omega_1)$.
Does $X \sub L^p(\mu)$ for some finite measure~$\mu$? 
\end{problem}

A Banach space~$X$ is said to be {\em Hilbert generated} if there exist a Hilbert space~$H$ and an
operator $T:H \to X$ having dense range (see e.g. \cite[Section~6.3]{fab-alt-JJ}). The basic examples 
of Hilbert generated spaces are: separable spaces, $c_0(\Gamma)$ (for any set~$\Gamma$) and $L^p(\mu)$ for any $1\leq p\leq 2$
and any finite measure~$\mu$ (in this case we can take $H=L^2(\mu)$ and $T$ the inclusion operator).
Every super-reflexive space is a {\em subspace} of a Hilbert generated space and there are examples
of super-reflexive spaces which are not Hilbert generated, \cite{fab-alt-J-4}. However, it seems to be unknown if every space
$L^p(\mu)$ with $2<p<\infty$ (and a finite measure~$\mu$) is Hilbert generated. Of course, this is true 
if $L^p(\mu)$ is separable, but also if ${\rm dens}(L^p(\mu))=\omega_1$, see~\cite{fab-alt-J-4}.

\begin{problem}\label{problem:LpHilbertGenerated}
Let $2<p<\infty$ and let $\mu$ be a finite measure. Is $L^p(\mu)$ Hilbert generated?
\end{problem}

\subsection{Compactness in Banach spaces}\label{subsection:compacidad}

A Banach space $X$ is said to be {\em weakly compactly generated} (WCG) if
there is a weakly compact set $K \sub X$ such that $X=\overline{{\rm span}}(K)$.
This class of Banach spaces includes all separable spaces, reflexive spaces, Hilbert generated spaces, etc. and
plays a very important role in non-separable Banach space theory (see e.g. \cite[Chapter~13]{fab-ultimo} and 
\cite[Chapter~6]{fab-alt-JJ}). In the particular case of a Banach {\em lattice}~$X$, Diestel~\cite{die6} asked 
whether the property of being WCG is equivalent to the existence of a weakly compact set 
$L \sub X$ such that the {\em sublattice} generated by~$L$ is dense in~$X$.
In~\cite{avi-alt-3} we answered in the affirmative 
this question for order continuous Banach lattices, but the general case still remains open.

\begin{problem}[Diestel, \cite{die6}]\label{problem:Diestel}
Let $X$ be a Banach lattice for which there is a weakly compact set $L \sub X$ such that
the sublattice generated by~$L$ is dense in~$X$. Is $X$ WCG?
\end{problem}

The concept of {\em weakly precompactly generated} (WPG) Banach space
was introduced by Haydon in~\cite{hay10}. A subset $C$ of a Banach space is called
{\em weakly precompact} if every sequence in~$C$ admits a weakly Cauchy subsequence or,
equivalently, if $C$ is bounded and contains no sequence which is 
equivalent to the usual basis of~$\ell^1$ (thanks to Rosenthal's $\ell^1$-theorem, \cite{ros-J-3}). 
A Banach space $X$ is said to be WPG if there is a weakly precompact set $C \sub X$ 
such that $X=\overline{{\rm span}}(C)$.
This class includes all WCG spaces and all Banach spaces not containing~$\ell^1$. 
Some results on these spaces can also be found in~\cite[Section~2.3]{sch-PhD}. 
While preparing his PhD Thesis, G.~Martínez-Cervantes~\cite{gon2} is studying WPG spaces and the 
class of {\em weak Radon-Nikod\'{y}m} (WRN) compacta introduced in~\cite{gla-meg-1}. 
A compact $K$ is said to be WRN if it is homeomorphic to a $w^*$-compact subset of the dual
of a Banach space not containing~$\ell^1$; this condition is equivalent to $C(K)$ being WPG.
The class of WRN compacta is larger than the well-known class of Radon-Nikod\'{y}m compacta. 
Recall that a compact $K$ is called {\em Radon-Nikod\'{y}m} if it is
homeomorphic to a $w^*$-compact subset of the dual of an Asplund space.

In~\cite{hay10} Haydon asked, using the language of WPG Banach spaces, whether every WRN (infinite) compact 
admits a convergent (non-stationary) sequence. This question is related to ``Efimov's problem'' on the existence of
(infinite) compacta without subsets homeomorphic to~$\beta\N$ and without convergent (non-stationary) sequences. 
Until now such examples have been constructed only under additional axioms of set theory (see e.g. \cite{dow-she,har}).

\begin{problem}[Haydon, \cite{hay10}]\label{problem:Haydon}
Does every WRN (infinite) compact admit a convergent (non-stationary) sequence?
\end{problem}

In~\cite{avi-ple-rod-5} we have studied the partially ordered set $\cK(B_X)$ 
of all weakly compact subsets of the closed unit ball of a separable Banach space~$X$
(the order being given by inclusion). To measure the complexity of~$\cK(B_X)$ we use Tukey ordering, which is a useful tool to isolate
essential properties of ordered structures in measure theory and topology, see \cite{fre12,fre13,sol-tod}. For instance,
Fremlin~\cite{fre12} proved that, given a coanalytic separable metric space~$E$, the partially ordered set 
$\cK(E)$ of all compact subsets of~$E$ is Tukey equivalent to one of the following:
$\{0\}$, $\N$, $\N^\N$, $\cK(\mathbb{Q})$. In~\cite{avi-ple-rod-5} we use advanced
techniques of descriptive set theory to show, under the axiom of analytic determinacy, 
that $\cK(B_X)$ is Tukey equivalent to $\{0\}$, $\N^\N$, $\cK(\mathbb{Q})$ or~$[\mathfrak{c}]^{<\mathbb{N}}$
(the family of all finite subsets of~$\mathfrak{c}$). This classification result is valid in ZFC
if $X\not\supseteq \ell^1$, but we do not know what happens in general.

\begin{problem}\label{problem:Tukey}
Is it relatively consistent that there is a separable Banach space~$X$ such that $\cK(B_X)$ is not
Tukey equivalent to $\{0\}$, $\N^\N$, $\cK(\mathbb{Q})$ or~$[\mathfrak{c}]^{<\mathbb{N}}$?
\end{problem}

\subsection{$w^*$-null sequences in dual spaces}\label{subsection:SequencesDual}

A subset $C$ of a Banach space~$X$ is said to be {\em limited} if $\lim_{n\to \infty} \, \sup_{x\in C}|x_n^*(x)|=0$ 
for every $w^*$-null sequence $(x_n^*)$ in~$X^*$. This concept is intimately related to different
notions of compactness. It is easy to check that every relatively norm compact set is limited. 
Bourgain and Diestel~\cite{bou-die} proved that, in general, every limited set is weakly precompact; they also showed that
every limited subset of a Banach space not containing~$\ell^1$ is relatively weakly compact. This last result extends to any WPG space, as
pointed out in~\cite[Cor.~2.3.3]{sch-PhD}. In another direction, 
a Banach space is said to have the {\em Gelfand-Phillips property} if every limited subset is relatively norm compact
(see e.g. \cite{bou-die,dre2,sch-2}). For instance, every Banach space having $w^*$-sequentially compact dual ball (e.g. 
WCG or Asplund spaces, see \cite[Chapter~XIII]{die-J})
or having the separable complementation property enjoys the Gelfand-Phillips property. Recall
that a Banach space $X$ has the {\em separable complementation property} (SCP) if
every separable subspace of~$X$ is contained in a separable complemented subspace of~$X$. 
Typical examples of Banach spaces with the SCP are WCG spaces and Banach lattices not containing~$c_0$. 
The following open problem is directly 
connected with a question on Pettis integration raised by Talagrand~\cite[4-2-6]{tal} (see Subsection~\ref{subsection:Pettis}):

\begin{problem}\label{problem:GelfandPhillips}
Let $X$ be a Banach space not containing $c_0$. Does $X$ have the Gelfand-Phillips property? 
\end{problem}

Mazur's theorem ensures that every weakly convergent sequence in a Banach space
admits a convex block subsequence which converges in norm. 
Recall that a {\em convex block subsequence} of a sequence $(x_n)$ in a Banach space is
a sequence $(y_k)$ of vectors of the form $y_k=\sum_{n\in I_k}a_n x_n$, where $I_1,I_2,\dots$ 
are finite subsets of~$\N$ with $\max(I_k)<\min(I_{k+1})$ and the scalars $a_n\geq 0$ satisfy 
$\sum_{n\in I_k}a_n=1$ for every $k\in \N$. 
A Banach space~$X$ has {\em property~(K)} if every $w^*$-null sequence $(x_n^*)$ in~$X^*$ admits
a convex block subsequence $(y^*_k)$ which converges in the Mackey topology $\mu(X^*,X)$, that is,
$\lim_{k\to \infty} \, \sup_{x\in K}|y_k^*(x)|=0$ for every weakly compact set $K \sub X$.
This concept, attributed to Kwapie\'{n}, was used in \cite{kal-pel} to study certain questions about 
subspaces of~$L^1[0,1]$. A variant of property~(K) was employed in~\cite{fig-alt2} to prove
that, in general, the SCP is not inherited by subspaces. The basic examples of Banach spaces with property~(K)
are the reflexive ones and, more generally, Grothendieck spaces (e.g. $\ell^\infty$) and the strongly WCG spaces of~\cite{sch-whe}
(e.g. $L^1(\mu)$ for any finite measure~$\mu$). A Banach space $X$ is said to be {\em strongly WCG} (SWCG) 
if there is a weakly compact set $K_0 \subseteq X$ such that, for every weakly compact $K \subseteq X$ and every 
$\varepsilon>0$, we have $K \subseteq nK_0 + \varepsilon B_X$ for some $n\in \N$. 

Pe\l czy\'{n}ski (see e.g.~\cite{fig-alt2,fra-ple}) proved that the $\ell^1$-sum of $\mathfrak{b}$ copies of~$L^1[0,1]$ fails
property~(K). Extending a result of~\cite{fra-ple}, we have recently shown that property~(K)
is preserved by $\ell^1$-sums of less than~$\mathfrak{p}$ summands, \cite{avi-ple-rod-6}.
For more information about cardinals $\mathfrak{p}$ and $\mathfrak{b}$ 
(that satisfy $\omega_1\leq \mathfrak{p}\leq \mathfrak{b}\leq \mathfrak{c}$), see e.g.~\cite{dou2}.

\begin{problem}\label{problem:l1sumasK}
What is the least cardinality of a family of Banach spaces with property~(K) whose $\ell^1$-sum fails such property?
\end{problem}

\begin{problem}\label{problem:BaseIncondicionalK}
Let $X$ be a separable Banach space which is weakly sequentially complete. Does $X$ have property~(K)?
\end{problem}

\section{Problems in measurability in Banach spaces}\label{section:measurability}

\subsection{Baire $\sigma$-algebras}\label{subsection:Baire}

The \emph{Baire $\sigma$-algebra} of a topological space~$T$, denoted by $\Ba(T)$, is the one generated by all continuous functions
from~$T$ to~$\erre$. It is contained in the \emph{Borel $\sigma$-algebra} of~$T$ (denoted by $\Bo(T)$) and, in general, such
inclusion is strict. If $E$ is a locally convex space, a result of Edgar~\cite{edgar1}
says that the Baire $\sigma$-algebra of~$(E,weak)$
is the $\sigma$-algebra generated by all elements of the topological dual of~$E$. In particular, 
if $X$ is a Banach space and $\Gamma \sub X^*$ is a total set, then
$\Ba(X,\sigma(X,\Gamma))$ is the $\sigma$-algebra generated by~$\Gamma$; of course, one has 
$\Ba(X,\sigma(X,\Gamma)) \sub \Ba(X,w)$. A Banach space~$X$ is said to have {\em property $\cD$} 
if $\Ba(X,\sigma(X,\Gamma)) = \Ba(X,w)$ for every total set~$\Gamma \sub X^*$.
This property is fulfilled by every Banach space~$X$ such that $(X^*,w^*)$ is {\em angelic}
(e.g. if $X$ is WCG),~\cite{gul-J}. The angelicity of $(X^*,w^*)$ means that every element of the $w^*$-closure of a bounded set $A\sub X^*$ is the $w^*$-limit of a sequence contained in~$A$. The space of Johnson-Lindenstrauss $JL_2$ has property~$\cD$
but fails to have $w^*$-angelic dual,~\cite{pli3}. A technical condition between the $w^*$-angelicity of the dual and property $\cD$ is the following: 
\begin{itemize}
\item[$\cD'$:] every $w^*$-sequentially closed linear subspace of~$X^*$ is $w^*$-closed.
\end{itemize} 

\begin{problem}[Plichko, \cite{pli3}]\label{problem:Plichko}
Is property $\cD'$ equivalent to property~$\cD$?
\end{problem}

Let $K$ be compact. The dual space $C(K)^*$ is identified with the space of all regular Borel signed measures
on~$K$. The subset of $C(K)^*$ consisting of all probabilities is denoted by~$P(K)$. 
Let $\mathfrak{T}_p$ be the topology on~$C(K)$ of pointwise convergence on~$K$.
In view of the above, the $\sigma$-algebra $\Ba(C(K),\mathfrak{T}_p)$ is generated by the set of 
evaluation functionals $\Delta_K:=\{\delta_t:t\in K\} \sub P(K)$, where 
$\delta_t(f):=f(t)$ for all $t\in K$ and $f\in C(K)$.
The equality $\Ba(C(K),\mathfrak{T}_p)=\Ba(C(K),w)$ holds true if every $\mu \in P(K)$
admits a {\em uniformly distributed sequence}, i.e. a sequence $(t_n)$ in~$K$ such that, for every $f\in C(K)$,
$$
	\int_K f \, d\mu=\lim_{m\to \infty}\frac{1}{m}\sum_{n=1}^m f(t_n).
$$
This happens in several cases, for instance, when $K=\{0,1\}^{\mathfrak{c}}$ or when
$K$ is Eberlein, Radon-Nikod\'{y}m, etc. (see~\cite[\S 491]{freMT-4} and \cite{mer-J}). Recall that 
$K$ is said to be {\em Eberlein} if it is homeomorphic to a weakly compact subset of a Banach space.

However, the coincidence of $\Ba(C(K),\mathfrak{T}_p)$ and $\Ba(C(K),w)$ does not guarantee the existence
of uniformly distributed sequences for every $\mu\in P(K)$, at least under the Continuum Hypothesis~\cite{avi-ple-rod-4}.
Now, write ${\rm co}(\Delta_K)$ for the convex hull
of~$\Delta_K$ and let ${\rm Seq}({\rm co}(\Delta_K))$ be its $w^*$-sequential closure in~$C(K)^*$
(i.e. the smallest $w^*$-sequentially closed subset of~$C(K)^*$ containing~${\rm co}(\Delta_K)$). 

\begin{problem}\label{problem:BaireSequential}
Let $K$ be compact such that $\Ba(C(K),\mathfrak{T}_p)=\Ba(C(K),w)$. Does the equality
${\rm Seq}({\rm co}(\Delta_K)) = P(K)$ hold?
\end{problem}

An equivalent norm $\|\cdot\|$ in a Banach space~$X$ is $\Ba(X,w)$-measurable (as a function from $X$ to~$\erre$) 
if and only if its balls belong to~$\Ba(X,w)$. The following implications are easy to check:
\begin{center}
$B_{X^*}$ is $w^*$-separable $\Longrightarrow$ $\|\cdot\|$ is $\Ba(X,w)$-measurable $\Longrightarrow$ 
$X^*$ is $w^*$-separable. 
\end{center}
The converse implications do not hold in general, as we shown in~\cite{rod9} by using certain 
equivalent norms in~$\ell^\infty$. However, we do not know what happens with the converse implications
when $X$ is a $C(K)$ space equipped with the {\em supremum norm}~$\|\cdot\|_\infty$. Under the
Continuum Hypothesis, Talagrand~\cite{tal12} constructed a compact~$K$
such that $C(K)^\ast$ is $w^\ast$-separable while $B_{C(K)^\ast}$ is not.
In~\cite{avi-ple-rod-3} we gave a ZFC example of a compact with these properties and we also proved
some partial results about the measurability of $\|\cdot\|_\infty$ in the corresponding $C(K)$ space.

\begin{problem}\label{problem:NormaSupremo}
Let $K$ be compact. Is the measurability of $\|\cdot\|_\infty$
with respect to $\Ba(C(K),w)$ equivalent to the $w^*$-separability of $B_{C(K)^*}$ or~$C(K)^*$?
\end{problem}

\subsection{Borel $\sigma$-algebras}\label{subsection:Borel}

Given a compact~$K$, the space $C(K)$ admits several topologies (of uniform, weak or pointwise convergence) 
leading to different $\sigma$-algebras:
\begin{equation}\label{eqn:sigma-algebras}
\begin{array}{c c c c c}
\Ba(C(K),\mathfrak{T}_p) & \subset & \Bo(C(K),\mathfrak{T}_p) & \mbox{ } & \mbox{ } \\
 \cap & \mbox{ } & \cap & \mbox{ } & \mbox{ } \\
\Ba(C(K),w) & \subset & \Bo(C(K),w) & \subset & \Bo(C(K))
\end{array}
\end{equation}
All these $\sigma$-algebras coincide if $K$ is metrizable, but also in some cases beyond metrizability, 
like $K=\{0,1\}^{\omega_1}$,~\cite{avi-ple-rod-2}.
On the other hand, $K=\beta\N$ is an example for which all the $\sigma$-algebras of~\eqref{eqn:sigma-algebras}
are different, see \cite{mar-ple,tal9}. In general, the equalities between such $\sigma$-algebras are 
closely related with topological properties of~$K$ and Banach space properties of~$C(K)$. 
For instance, it is known that the equality $\Bo(C(K),\mathfrak{T}_p) = \Bo(C(K))$ holds true if $K$ is an Eberlein compact or, more generally,
a Valdivia compact (see e.g. \cite{edgar1,edgar2} and~\cite[Chapter~VII]{dev-alt-J}). 
The following questions remain open:

\begin{problem}[Burke-Pol, \cite{bur-pol-1}]\label{problem:Burke-Pol}
Let $K$ be compact such that 
$$
	\Bo(C(K),w) \neq \Bo(C(K)).
$$ 
Is there a norm discrete subset of~$C(K)$ which does not belong to $\Bo(C(K),w)$?
\end{problem}

\begin{problem}[Marciszewski-Pol, \cite{mar-pol-2}]\label{problem:Marciszewski-Pol}
Let $K$ be compact such that 
$$
	\Bo(C(K),w) = \Bo(C(K)).
$$ 
Does the equality $\Bo(C(K),\mathfrak{T}_p)=\Bo(C(K))$ hold?
\end{problem}

A topological space $T$ is said to be {\em Radon} if every probability $\mu$ on~$\Bo(T)$ 
is Radon (i.e. $\mu(E)=\sup\{\mu(L): \, L \sub E, \, L \mbox{ compact}\}$ for every $E\in \Bo(T)$). 
A classical result states that a complete metric space~$T$ is Radon if ${\rm dens}(T)$ 
is a cardinal of measure zero (e.g. if ${\rm dens}(T)\leq \omega_1$). L.~Schwartz asked if $(X,w)$ is Radon for any
Banach space~$X$. If $X$ is separable, then $\Bo(X,w)=\Bo(X)$ and therefore $(X,w)$ is Radon;
the same happens if $X$ is WCG and ${\rm dens}(X)=\omega_1$.
However, there exist Banach spaces with density character~$\mathfrak{c}$ which are not Radon with its weak topology, like $\ell^\infty/c_0$. 
The references \cite[\S 466]{freMT-4}, \cite{jay-alt10} and~\cite{tal} contain more information about this topic. 
For any Banach space~$X$, every Radon probability on~$\Bo(X,w)$ extends to a Radon probability on~$\Bo(X)$
(Phillips-Grothendieck). Thanks to this result,
an affirmative answer to the following question would imply that $(\ell^\infty,w)$ is not Radon,
thus solving a long-standing conjecture:

\begin{problem}[Fremlin, \cite{freMT-4}]\label{problem:FremlinRadon}
Is there a probability on $\Bo(\ell^\infty,w)$ which cannot be extended to a probability on $\Bo(\ell^\infty)$?
\end{problem}

\subsection{Measurable selectors}\label{subsection:selectors}

Let $(\Omega,\Sigma,\mu)$ be a complete probability space and $X$ a Banach space.
Denote by $\mathcal{P}_0(X)$ the family of all non-empty subsets of~$X$.
A \emph{multi-function} is a map $F: \Omega \to \mathcal{P}_0(X)$ and 
a \emph{selector} of~$F$ is a function $f: \Omega \to X$ such that $f(t) \in F(t)$ for every $t \in \Omega$. 
Most of the results about the existence of ``measurable'' selectors
are restricted to the context of separable spaces, that allows to use 
classical descriptive set theory. For instance, the Kuratowski
and Ryll-Nardzewski selection theorem~\cite{kur-ryl} ensures that, for separable~$X$, a multi-function $F:\Omega \to \mathcal{P}_0(X)$ 
has $\Bo(X)$-measurable selectors if it takes closed values and satisfies:
$$
	F^{-}(G):=\{t\in \Omega: \, F(t)\cap G \neq \emptyset\}\in \Sigma\quad
	\mbox{for every open }G \sub X.
$$
A multi-function $F:\Omega \to \mathcal{P}_0(X)$ is said to be {\em scalarly measurable} if $F^{-}(H)\in \Sigma$ for every 
open half-space $H\sub X$. For arbitrary Banach spaces, 
we showed in~\cite{cas-kad-rod-3} that every scalarly measurable multi-function with {\em weakly compact} values 
admits scalarly measurable selectors (cf. \cite{fre-CKR}). Recall that 
a function $f:\Omega\to X$ is said to be {\em scalarly measurable} if the composition $x^*f:\Omega \to \erre$ 
is measurable for every $x^*\in X^*$ (equivalently, $f$ is $\Ba(X,w)$-measurable).
The result of~\cite{cas-kad-rod-3} allowed to develop the set-valued Pettis integral theory
in non-separable Banach spaces, \cite{cas-kad-rod-2,mus9,mus11}. 

\begin{problem}\label{problem:multi}
Let $F:\Omega \to \mathcal{P}_0(X)$ be a scalarly measurable multi-function with convex closed values.
Does $F$ admit scalarly measurable selectors?
\end{problem}

A {\em set selector} of a family of sets~$\cC \sub \mathcal{P}_0(X)$ is a 
map $\psi:\cC \to X$ such that $\psi(C)\in C$ for every $C\in \cC$.
If $F:\Omega \to \cC$ is any multi-function, then
the composition $\psi\circ F:\Omega \to X$ is a selector of~$F$, and it is natural to ask whether
such selector enjoys some nice measurability or continuity property provided $F$ does. Some important theorems 
about measurable selectors (Kuratowski and Ryll-Nardzewski) or first Baire class selectors
(Jayne and Rogers) can be deduced from the existence of suitable set selectors,~\cite{gho-alt2}. 
We say that the family $\cC$ has {\em property~(SC)} if there is
a set selector $\psi:\cC\to X$ such that $\psi\circ F$ is scalarly measurable
for every scalarly measurable multi-function $F:\Omega \to \cC$. 
Denote by $wk(X)$ the family of all weakly compact 
non-empty subsets of~$X$. Then $wk(X)$ has property~(SC) if $X^*$ is $w^*$-separable 
or if $(X^*,w^*)$ is angelic and has density character~$\omega_1$, 
see \cite{vld} and~\cite{cas-kad-rod-2} (cf. \cite{fre-CKR}).
Fremlin~\cite{fre-CKR} proved that, in general, $wk(X)$ may fail property (SC), even for Hilbert spaces. 
The following questions remain open: 

\begin{problem}[Fremlin, \cite{fre-CKR}]\label{problem:selectorsets1}
Suppose ${\rm dens}(X)=\omega_1$. Does $wk(X)$ have property (SC)?
\end{problem}

\begin{problem}[Fremlin, \cite{fre-CKR}]\label{problem:selectorsets2}
Does the family $cwk(X)$ of all convex weakly compact 
non-empty subsets of~$X$ have property (SC)?
\end{problem}

\section{Problems in vector integration}\label{section:integrals}

\subsection{Riemann integral}\label{section:Riemann}

Lebesgue's criterion of Riemann integrability does not work in general for 
vector-valued functions: there exist Riemann integrable functions 
defined on~$[0,1]$ with values in a Banach space which are not continuous a.e. (see e.g. \cite{gor-2}). 
This phenomenon occurs in almost all classical Banach spaces. A Banach space~$X$ is said to have 
the {\em Lebesgue property} (LP) if every Riemann integrable function $f:[0,1] \to X$ is continuous a.e. 
For instance, $\ell^{1}(\Gamma)$ has the LP for any set~$\Gamma$. A complete characterization of Banach spaces having the~LP
is missing, although there are partial results related to fine properties of the structure of Banach spaces.
A result attributed to da Rocha and Pe\l czy\'{n}ski asserts that every asymptotic $\ell^1$ Banach space 
(e.g. Tsirelson's space) has the LP. On the other hand, Haydon~\cite{hay9} proved that 
Banach spaces having the LP share a feature with Banach spaces having the Schur property, namely:
every spreading model is equivalent to the usual basis of~$\ell^1$. 
However, there exist separable Banach spaces having the Schur property and failing the LP,~\cite{hay9}.
The paper~\cite{nar} contains alternative proofs of these results.

\begin{problem}\label{problem:LP}
Characterize Banach spaces having the LP.
\end{problem}

Alexiewicz and Orlicz~\cite{ale-orl-J} showed an example of a Riemann integrable function with values in~$C[0,1]$
without points of continuity with respect to the weak topology. A Banach space~$X$ is said to have the {\em weak Lebesgue property} (WLP)
if every Riemann integrable function $f:[0,1] \to X$ is {\em weakly} continuous a.e. 
It is easy to see that all spaces with separable dual have the WLP. Beyond this case, the space 
$L^{1}[0,1]$ has the WLP,~\cite{wan-2}. This result was extended to certain Lebesgue-Bochner spaces in~\cite{cal-alt-3} 
and, recently, it has been improved by G.~Martínez-Cervantes~\cite{gon}, showing that $C(K)^*$ has the WLP for a wide class 
of compacta $K$ including those which are Eberlein or Radon-Nikod\'{y}m. Some partial answers to the following question
have also been given in~\cite{gon}:

\begin{problem}\label{problem:WLP}
For which cardinals~$\kappa$ does the space $\ell^2(\kappa)$ have the WLP?
\end{problem}

Given a Banach space~$X$ and a bounded function $f: [0,1] \to X$, non-necessarily Riemann integrable,
one can consider the subset $I(f)$ of~$X$ consisting of all limits 
of sequences of Riemann sums of~$f$ based on partitions with diameter tending to~$0$.
The properties of the set $I(f)$ have been widely studied by several authors, see \cite[Appendix]{kad-kad-J}. 
It is known that $I(f)\neq \emptyset$ if $X$ is separable or if $X$ is super-reflexive; in the second case, 
the set $I(f)$ is convex. Recently a similar concept of
``limit set'' based on the Henstock-Kurzweil integral has been discussed in~\cite{cap-alt}. One of
the open questions along this line is:

\begin{problem}[Kadets-Kadets, \cite{kad-kad-J}]\label{problem:KadetsKadets}
Is the set $I(f)$ convex for any bounded function $f:[0,1]\to c_0$?
\end{problem}

\subsection{Pettis integral}\label{subsection:Pettis}

Let $(\Omega,\Sigma,\mu)$ be a complete probability space and $X$ a Banach space. 
A function $f:\Omega \to X$ is said to be {\em scalarly bounded} 
if there is $M>0$ such that for each $x^*\in B_{X^*}$ we have $|x^*f|\leq M$ $\mu$-a.e. 
A function $f:\Omega \to X$ is called {\em Pettis integrable} if: (i)~$x^{*}f$ 
is integrable for every $x^{*}\in X^{*}$; (ii)~for each $A\in \Sigma$
there is a vector $\int_A f \, d\mu \in X$ such that $x^{*}(\int_A f \, d\mu)=\int_{A}x^{*} f \, d\mu$
for all $x^{*}\in X^{*}$. The Pettis integral theory is closely related to deep questions about 
pointwise compact sets of measurable real-valued functions (see e.g. \cite{tal}).
The space~$X$ has the {\em Pettis integrability property} (PIP)
if every scalarly bounded and scalarly measurable function taking values in~$X$ is Pettis integrable.
For instance, if $(X^*,w^*)$ is angelic, then $X$ has the PIP. On the other hand,
Fremlin and Talagrand~\cite{fre-tal} proved that $\ell^\infty$ fails the PIP. 

For any Pettis integrable function $f:\Omega \to X$, the map $\nu_f:\Sigma \to X$
defined by $\nu_f(A):=\int_A f \, d\mu$ is countably additive and, in particular, its range
$$
	\mathcal{R}(\nu_f)=\Big\{\int_A f \, d\mu: \ A\in \Sigma\Big\} \sub X
$$
is relatively weakly compact (see Subsection~\ref{subsection:range}). 
One can approximate $f$ nicely by simple functions whenever the
set~$\mathcal{R}(\nu_f)$ is separable or, even better, relatively norm compact.
The paper~\cite{fre-tal} exhibited the first example of a Pettis integrable function~$f$ such that 
$\mathcal{R}(\nu_f)$ is not relatively norm compact. Such examples are very pathological, 
since the relative norm compactness of~$\mathcal{R}(\nu_f)$ is guaranteed in many cases, for instance:
if $\mu$ is a perfect measure (e.g. a Radon measure), if $X\not\supseteq \ell^{1}(\omega_{1})$
(e.g. if $(X^*,w^*)$ is angelic), if $X$ is the dual of a Banach space not containing~$\ell^1$, etc. (see e.g. \cite{Musial,mus3,tal}).
There are no characterizations of those Banach spaces~$X$ and complete probability spaces~$(\Omega,\Sigma,\mu)$ for which
$\mathcal{R}(\nu_f)$ is relatively norm compact (resp. separable) for every Pettis integrable function $f:\Omega \to X$.

\begin{problem}[Plebanek, \cite{ple}]\label{problem:Plebanek}
Suppose $X$ has the PIP. Is $\mathcal{R}(\nu_f)$ separable for every Pettis integrable function $f:\Omega \to X$?
\end{problem}

A \emph{quasi-Radon} probability space (see e.g. \cite[Chapter~41]{freMT-4}) is a quadruple
$(\Omega,\mathfrak{T},\Sigma,\mu)$, where $(\Omega,\Sigma,\mu)$ is a
complete probability space and $\mathfrak{T}\subset \Sigma$ is a topology
on~$\Omega$ such that $\mu$ is inner regular with respect to the
collection of all closed sets, and $\mu(\bigcup \mathcal{G})=\sup\{\mu(G):
G\in \mathcal{G}\}$ for every upwards directed family $\mathcal{G}\subset
\mathfrak{T}$. Of course, every Radon probability space is quasi-Radon.

\begin{problem}[Fremlin, \cite{fre5}]\label{problem:FremlinQuasiRadon}
Suppose $\mu$ is quasi-Radon. Is $\mathcal{R}(\nu_f)$ relatively norm compact for
every Pettis integrable function $f:\Omega \to X$?
\end{problem}

It is easy to check that all sets of the form $\mathcal{R}(\nu_f)$ are limited. Therefore, 
an affirmative answer to Problem~\ref{problem:GelfandPhillips} 
(see Subsection~\ref{subsection:SequencesDual}) would solve the following:

\begin{problem}[Talagrand, \cite{tal}]\label{problem:Talagrand}
Suppose $X\not\supseteq c_0$. Is $\mathcal{R}(\nu_f)$ relatively norm compact for
every Pettis integrable function $f:\Omega \to X$?
\end{problem}

\subsection{McShane integral}\label{subsection:McShane}

Let $X$ be a Banach space. A function $f:[0,1] \to X$ is said to be {\em McShane integrable},
with integral $x\in X$, if it satisfies the following condition: for every $\epsilon>0$ there is 
$\delta:[0,1] \to (0,\infty)$ such that
$\|\sum_{j=1}^{p}\lambda(I_{j})f(t_{j})-x\|< \epsilon$ for every 
partition $\{I_1,\dots,I_p\}$ of~$[0,1]$ into
finitely many non-overlapping closed intervals and every choice of points $t_j\in [0,1]$ 
such that $I_j\sub (t_j-\delta(t_j),t_j+\delta(t_j))$.
The notion of McShane integrability lies between Bochner and Pettis integrability; 
moreover, McShane and Pettis integrability coincide for functions with values in separable Banach spaces,~\cite{fre-men}. 
Let $\mathcal{U}$ be the class of those Banach spaces~$X$ such that every Pettis integrable function $f:[0,1]\to X$ is McShane integrable. 
Within the context of non-separable Banach spaces, Di Piazza and Preiss~\cite{dip-pre} proved that 
super-reflexive spaces and $c_0(\Gamma)$ (for any set $\Gamma$) belong to~$\mathcal{U}$, raising the problem of
whether the same holds for any WCG space. In~\cite{rod10} we proved that the space $L^1(\mu)$ (for any 
finite measure~$\mu$) belongs to~$\mathcal{U}$. Later we showed that all (subspaces of) 
Hilbert generated spaces belong to the class~$\mathcal{U}$ (see~\cite{dev-rod}), thus generalizing all previous results in this direction. 
In general, the answer to the question of Di Piazza and Preiss is negative: in~\cite{avi-alt} we constructed Pettis integrable functions, with
values in reflexive Banach spaces, that are not McShane integrable. The solution uses techniques of infinite
combinatorics related to Fremlin's ``problem DU''~\cite{freDU}. 

Let~$\cK$ denote the class of all compacta~$K$ such that $C(K)\in \mathcal{U}$. The results of~\cite{avi-alt,dev-rod} show that
$\cK$ contains every uniform Eberlein compact, but not every Eberlein compact.
Recall that a compact~$K$ is said to be {\em uniform Eberlein} if
it is homeomorphic to a weakly compact subset of a Hilbert space or, equivalently, if $C(K)$ is Hilbert generated.
The paper~\cite{fab-mcs} provides another approach to the results of~\cite{avi-alt,dev-rod} and discusses
some classical examples of Eberlein compacta which are not uniform Eberlein, in connection with the following open problem:

\begin{problem}\label{problem:McShane}
Are there non-uniform Eberlein compacta in the class~$\cK$?
\end{problem}

The McShane integral can be generalized to the case of functions defined on a quasi-Radon 
probability space~\cite{fre5}. In this context the following question arises (it has affirmative answer
if $\mathfrak{T}_1$ has a countable basis~\cite{avi-alt}):

\begin{problem}[Fremlin, \cite{fre5}]\label{problem:FremlinMcShane}
Let $f:\Omega \to X$ be a function defined in a complete probability space $(\Omega,\Sigma,\mu)$ and
taking values in a Banach space~$X$. Let $\mathfrak{T}_1$ and $\mathfrak{T}_2$ be two topologies on~$\Omega$ for which
$\mu$ is quasi-Radon. Suppose $f$ is McShane integrable with respect to~$\mathfrak{T}_1$. Is $f$ McShane integrable
with respect to~$\mathfrak{T}_2$?
\end{problem}

\section{Problems in vector measures}\label{section:MedidasVectoriales}

\subsection{Range of a vector measure}\label{subsection:range}

By a {\em vector measure} we mean a countably additive map
$\nu:\Sigma\to X$ defined on a $\sigma$-algebra~$\Sigma$ and taking values in a Banach space~$X$. 
We denote by $|\nu|$ (resp. $\|\nu\|$) the variation (resp. semivariation) of~$\nu$.
A theorem of Bartle, Dunford and Schwartz~\cite{bar-alt} ensures that the range of~$\nu$, i.e. 
the set $\mathcal{R}(\nu)=\{\nu(E):E\in \Sigma\}$, is relatively weakly compact.
In fact, the range has the Banach-Saks property (every sequence in it admits a 
subsequence whose arithmetic means are norm convergent),~\cite{die-sei}.
It is natural to ask which subsets of a Banach space
(e.g. sequences with some convergence or summability property) are contained
in the range of a vector measure (see e.g.~\cite{sof}).
It is known that the closed unit ball $B_X$ of 
a Banach space~$X$ {\em is} the range of a vector measure if and only if
$X^*$ is isometric to a reflexive subspace of~$L^1(\mu)$ for some finite measure~$\mu$
(e.g. $X=L^p[0,1]$ or $X=\ell^p$ with $2\leq p < \infty$),~\cite{ana-die}.
However, it seems that there is no characterization of those Banach spaces~$X$ for which $B_X$ is {\em contained}
in the range of an $X$-valued vector measure. This class of spaces
is strictly contained in that of super-reflexive spaces,~\cite{cas-san}.

\begin{problem}\label{problem:ball}
Characterize the Banach spaces~$X$ for which $B_X$ is contained in the range of an $X$-valued vector measure.
\end{problem}

\begin{problem}[Anantharaman-Diestel, \cite{ana-die}]\label{problem:Sofi}
Let $K$ be a weakly compact subset of~$c_0$. Is $K$ contained in the range of a $c_0$-valued vector measure?
\end{problem}

Certain properties of a vector measure are determined by its range. For instance,
if $\nu_1$ and $\nu_2$ are vector measures such that ${\mathcal{R}}(\nu_1)={\mathcal{R}}(\nu_2)$, 
then $|\nu_1|$ is finite if and only if $|\nu_2|$ is finite,~\cite{rodpiazza-2}. 
In such case, $\nu_1$~is Bochner differentiable with respect to~$|\nu_1|$
if and only if $\nu_2$ is Bochner differentiable with respect to~$|\nu_2|$, \cite{rodpiazza-1}. 
An example of Fremlin and Talagrand~\cite{fre-tal} was used in~\cite{rod-rom}
to show that the previous result is not true in general for the Pettis integral. 

\begin{problem}[Rodríguez-Piazza - Romero-Moreno, \cite{rod-rom}]\label{problem:RPRM}
Let $\nu_1$ and $\nu_2$ be vector measures of finite variation such that ${\mathcal{R}}(\nu_1)={\mathcal{R}}(\nu_2)$.
Suppose further that $|\nu_1|$ and $|\nu_2|$ are perfect measures. Do $\nu_1$ and $\nu_2$ 
have the same Pettis differentiability character with respect to their variations?
\end{problem}

\subsection{$L^1$ spaces of a vector measure}\label{subsection:L1m}

Let $(\Omega,\Sigma)$ be a measurable space, $X$ a Banach space and
$\nu: \Sigma \to X$ a vector measure. A measurable function
$f:\Omega \to \erre$ is said to be {\em $\nu$-integrable} if: 
(i)~it is integrable with respect to the composition~$x^{*}\nu:\Sigma \to \erre$ for every $x^{*}\in X^{*}$; 
(ii)~for each $A \in \Sigma$ there is $\int_{A}f \, d\nu \in X$ such that
$x^{*}(\int_{A}f \, d\nu)=\int_{A}f \, d(x^{*}\nu)$ for all $x^{*}\in X^{*}$.
By identifying functions which coincide $\|\nu\|$-a.e. we obtain the Banach lattice $L^1(\nu)$
of all (equivalence classes of) $\nu$-integrable functions,
with the norm $\|f\|_{L^{1}(\nu)}=\sup\{\int_{\Omega}|f| \, d|x^{*}\nu|: x^{*} \in B_{X^{*}}\}$.
Every $L^1$ space of a vector measure is both WCG and a subspace of a Hilbert generated space, \cite{cur1,rod15}.
Several classical spaces can be represented in this way, since
any order continuous Banach lattice with a weak order unit is lattice isometric to the
$L^1$~space of a vector measure,~\cite{cur1,dep-alt}. 

In this theory the integration operator
$I_{\nu}:L^1(\nu) \to X$, $I_{\nu}(f):=\int_\Omega f \, d\nu$, plays an important role (see e.g. \cite{oka-alt}). 
In~\cite{rod15} we improved results of \cite{oka-alt2,oka-alt3} by showing that 
if $I_\nu$ is completely continuous and Asplund (e.g. if $I_\nu$ is compact or absolutely $p$-summing, $1\leq p <\infty$), 
then $|\nu|$ is finite and $L^1(\nu)=L^1(|\nu|)$; therefore, in such case the space $L^1(\nu)$
is an ``ordinary'' $L^1$~space. Recall that an operator between Banach spaces is said to be {\em Asplund} if it factors
through an Asplund space (for instance, every weakly compact operator is Asplund).
The particular case in which $X$ is an Asplund space
was considered in~\cite{cal-alt-5} to give a partial answer to the following question:

\begin{problem}[Okada - Ricker - Rodríguez-Piazza, \cite{oka-alt3}]
Let $X$ be a Banach space not containing~$\ell^1$ and let $\nu$ be an $X$-valued vector measure 
such that $I_\nu$ is completely continuous. Does $\nu$ have finite variation?
\end{problem}

Curbera~\cite{cur2} (cf. \cite{lip}) proved that if $\nu$ is non-atomic and $L^1(\nu)$ is separable, then
there is a vector measure $\tilde{\nu}:\Sigma\to c_0$ such that $L^1(\nu)=L^1(\tilde{\nu})$ with equal norms.
Lipecki~\cite{lip2} asked if there is a ``non-separable version'' of that result by using $c_0(\Gamma)$ as range space
for a suitable set~$\Gamma$. In~\cite{rod16} we have shown that if $\nu$ is non-atomic and ${\rm dens}(L^1(\nu))=\omega_1$, then
there is a vector measure $\tilde{\nu}:\Sigma\to \ell^\infty_c(\omega_1)$ such that $L^1(\nu)=L^1(\tilde{\nu})$ with equivalent norms.
Here $\ell^\infty_c(\omega_1)$ denotes the subspace of $\ell^\infty(\omega_1)$ made up of all countably supported vectors.

\begin{problem}[Lipecki, \cite{lip2}]\label{problem:Lipecki}
Suppose $\nu$ is non-atomic. Are there a set $\Gamma$ and a vector measure $\tilde{\nu}:\Sigma \to c_0(\Gamma)$
such that $L^1(\nu)=L^1(\tilde{\nu})$?
\end{problem}

For any $g \in B_{L^{\infty}(\nu)}$ and $x^{*}\in B_{X^{*}}$ we have a functional
$\gamma_{g,x^{*}} \in B_{L^{1}(\nu)^{*}}$ defined by the formula $\gamma_{g,x^{*}}(f):=\int_{\Omega} f g \, d(x^{*}\nu)$.
The set $\Gamma$ consisting of all such functionals
is norming and the topology $\sigma(L^1(\nu),\Gamma)$
is coarser than the weak topology of~$L^1(\nu)$. However, there are non-trivial connections
between both topologies, specially when $\mathcal{R}(\nu)$ is relatively norm compact. 
In such case, every $\sigma(L^1(\nu),\Gamma)$-convergent sequence is weakly convergent (see~\cite{lew4,oka})
and, if in addition $X$ is Asplund, then both topologies coincide on bounded sets~\cite{gra-rue}.
These statements can fail in the absence of some of the hypotheses, \cite{cur2,cur3}.
Along this line, it is relevant to know conditions ensuring that~$\Gamma$ is a {\em James boundary} of~$L^1(\nu)$, meaning that 
for every $f\in L^1(\nu)$ there is $\gamma_{g,x^{*}}\in \Gamma$ such that $\|f\|_{L^{1}(\nu)}=\gamma_{g,x^{*}}(f)$. 
This happens when $\mathcal{R}(\nu)$ is relatively norm compact or
$\nu$ is a positive vector measure with values in a Banach lattice,~\cite{cal-alt-6,man-J}. 
We stress that, in an arbitrary Banach space, every
bounded subset which is compact with respect to the topology of pointwise convergence on a James boundary 
is weakly compact: this striking result of Pfitzner~\cite{pfi-J} answered a long-standing question 
known as ``the boundary problem''. It is also worth mentioning that James boundaries
are useful to study summability in Banach spaces, see \cite{cal-alt-7,fon2-J,gas}.

\begin{problem}\label{problem:frontera}
Characterize when $\Gamma$ is a James boundary of $L^1(\nu)$.
\end{problem}

\section{Problems in Lebesgue-Bochner spaces}\label{section:Bochner}

Let $(\Omega,\Sigma,\mu)$ be a probability space, $X$ a Banach space and $1\leq p<\infty$. Denote by
$L^p(\mu,X)$ the Banach space of all (equivalence classes of) 
strongly measurable functions $f:\Omega\to X$ such that 
$$
	\|f\|_{L^p(\mu,X)}=\Big(\int_\Omega \|f(\cdot)\|^p \, d\mu\Big)^{1/p}<\infty.
$$
These spaces are usually called Lebesgue-Bochner spaces. We refer the reader to \cite{cem-men,die-uhl-J,lin-J}
for complete information on this topic.

\subsection{Topological properties}\label{subsection:TopologiaBochner}

Two important contributions on the weak topology of Lebesgue-Bochner spaces
are the characterization of weakly compact sets of Diestel, Ruess and Schachermayer~\cite{die-alt2}, 
and the parametric version of Rosenthal's $\ell^1$-theorem due to Talagrand~\cite{tal11}. 
Nowadays the topic still offers challenging problems. A generic question is as follows: {\em if $X$ satisfies
a certain property~(P), does $L^p(\mu,X)$ enjoy~(P) as well?}
The factorization theorem of Davis, Figiel, Johnson and Pe\l czy\'{n}ski paves the way to show
that $L^p(\mu,X)$ is WCG if $X$ is,~\cite{die7}. 
Talagrand~\cite{tal1} proved that the properties of being weakly K-analytic (WKA) and weakly countably determined (WCD)
also pass from~$X$ to $L^p(\mu,X)$. It seems to be unknown if the same happens
with weaker properties like having $w^*$-sequentially compact dual ball.

\begin{problem}\label{problem:ballSC}
Is $B_{L^p(\mu,X)^*}$ $w^*$-sequentially compact if $B_{X^*}$ is? 
\end{problem}

Schl\"{u}chtermann and Wheeler~\cite{sch-whe} asked if the space $L^1(\mu,X)$ is SWCG when $X$ is SWCG. The answer is affirmative for 
reflexive spaces and $L^1(\mu)$ (for any finite measure~$\mu$), but the general case remains open. 
In~\cite{laj-rod} we gave a positive answer to the problem when restricted
to weakly compact decomposable sets (e.g. sets of selectors of multi-functions).
Recently we provided related results by showing, for instance, that
for $1<p<\infty$ (and infinite-dimensional $L^p(\mu)$), the space $L^p(\mu,X)$ is a subspace of an SWCG space
if and only if $X$ is reflexive,~\cite{rod13}.
On the other hand, the ``strong'' versions of properties WKA and WCD were discussed in~\cite{kam-mer,kam-mer2,mer-sta-2}
and it is natural to ask if they pass to the corresponding Lebesgue-Bochner space.

\begin{problem}[Schl\"{u}chtermann-Wheeler, \cite{sch-whe}]\label{problem:SW}
Is the space $L^1(\mu,X)$ SWCG whenever $X$ is?
\end{problem}

\begin{problem}\label{problem:SWKA}
Is the space $L^p(\mu,X)$ strongy WKA (resp. strongly WCD) whenever $X$ is?
\end{problem}

\subsection{Structural properties}\label{subsection:EstructuraBochner}

Kwapie\'{n} proved that $L^p(\mu,X) \not\supseteq c_0$ if (and only if) $X\not \supseteq c_0$. 
For $1<p<\infty$ we have an analogous result about copies of~$\ell^1$, due to Bourgain, Maurey and Pisier. 
Nice proofs of these results can be found e.g. in \cite[Chapter~2]{cem-men}.
In another direction, Diestel and Uhl~\cite{die-uhl-J} asked if the space $L^p([0,1],X)$ has unconditional basis
when $1<p<\infty$ and $X$ has unconditional basis. Aldous~\cite{ald} answered in the negative this question by showing that
the existence of an unconditional basis in $L^p([0,1],X)$ implies that $X$ is super-reflexive. 

\begin{problem}\label{problem:Aldous}
Characterize those Banach spaces~$X$ such that $L^p([0,1],X)$ has unconditional basis ($1<p<\infty$).
\end{problem}

While $\ell^2$ and $L^2[0,1]$ are isomorphic, the previous result of Aldous implies that
the spaces $\ell^2(X)$ and $L^2([0,1],X)$ are not isomorphic in general (take $X=c_0$, for instance). On the other hand,
$\ell^2$ is isomorphic to a subspace of~$L^1[0,1]$, but in general $\ell^2(X)$ is not isomorphic to a
subspace of~$L^1([0,1],X)$. An example of such situation is obtained by taking $X=\ell^1$, since
$L^1([0,1],\ell^1)$ is SWCG but $\ell^2(\ell^1)$ is not a subspace of an SWCG space (see the previous subsection).

\begin{problem}[Diestel, \cite{die5}]\label{problem:DiestelQ7}
For which separable Banach spaces~$X$ the spaces $\ell^2(X)$ and $L^2([0,1],X)$ are isomorphic?
\end{problem}

\begin{problem}[Diestel, \cite{die5}]\label{problem:DiestelQ8}
For which separable Banach spaces~$X$ the space $L^1([0,1],X)$ 
contains a copy of $\ell^2(X)$?
\end{problem}

The Banach space $L^\infty(\mu,X)$, consisting of all (equivalence classes of) essentially bounded and strongly measurable functions,
also plays a relevant role in the theory (see e.g. \cite[Chapter~5]{cem-men}). While $\ell^\infty$ and $L^\infty[0,1]$ are isomorphic,
in general the spaces $\ell^\infty(X)$ and $L^\infty([0,1],X)$ are not. For instance,
$L^\infty([0,1],\ell^\infty)$ contains a {\em complemented} copy of~$c_0$, 
so it cannot be isomorphic to~$\ell^\infty(\ell^\infty) \simeq \ell^\infty$.

\begin{problem}[Cembranos-Mendoza, \cite{cem-men}]\label{problem:CM}
When are $\ell^\infty(X)$ and $L^\infty([0,1],X)$ isomorphic?
\end{problem}

\subsection{Scalar convergence}\label{subsection:DilworthGirardi}

Edgar~\cite{edg6} proved that $X$ is Asplund if and only if the following condition holds:
every sequence of strongly measurable functions $f_n:[0,1] \to X$ which is bounded 
in $L^\infty([0,1],X)$ and {\em converges scalarly} to~$0$ a.e. 
(i.e. for every $x^* \in X^*$ the sequence $(x^*f_n)$ converges to~$0$ a.e.)
is weakly convergent to~$0$ a.e., that is, there is a $\lambda$-null set $E\sub [0,1]$ such that
the sequence $(f_n(t))$ is weakly null in~$X$ for every $t\in [0,1]\setminus E$. On the other hand,
Dilworth and Girardi~\cite{dil-gir} studied the scalar version of the classical result relating
convergence in measure with a.e. convergence. Namely, they discussed the following property (for $1\leq p \leq \infty$):
\begin{itemize}
\item[($D_p$)] Every sequence of strongly measurable functions $f_n:[0,1] \to X$ which is bounded in $L^p([0,1],X)$ 
and converges {\em scalarly in measure} to~$0$
(i.e. for every $x^* \in X^*$ the sequence $(x^*f_n)$ converges to~$0$ in measure) admits a subsequence that converges 
scalarly to~$0$ a.e.
\end{itemize}
A couple of open questions along this line are:

\begin{problem}[Dilworth-Girardi, \cite{dil-gir}]\label{problem:DG1}
Does $L^1[0,1]$ have property ($D_\infty$)?
\end{problem}

\begin{problem}[Dilworth-Girardi, \cite{dil-gir}]\label{problem:DG2}
Are there reflexive (infinite-dimensional) Banach spaces with property ($D_1$)?
\end{problem}

\subsection*{Acknowledgements}
Research partially supported by {\em Ministerio de Econom\'{i}a y Competitividad - FEDER} (project MTM2014-54182-P).
This work was also partially supported by the research project 19275/PI/14 funded by {\em Fundaci\'{o}n S\'{e}neca - Agencia de Ciencia y Tecnolog\'{i}a 
de la Regi\'{o}n de Murcia} within the framework of {\em PCTIRM 2011-2014}.

\def\cprime{$'$}\def\cdprime{$''$}
  \def\polhk#1{\setbox0=\hbox{#1}{\ooalign{\hidewidth
  \lower1.5ex\hbox{`}\hidewidth\crcr\unhbox0}}} \def\cprime{$'$}
\providecommand{\bysame}{\leavevmode\hbox to3em{\hrulefill}\thinspace}
\providecommand{\MR}{\relax\ifhmode\unskip\space\fi MR }
\providecommand{\MRhref}[2]{%
  \href{http://www.ams.org/mathscinet-getitem?mr=#1}{#2}
}
\providecommand{\href}[2]{#2}

\bibliographystyle{amsplain}

\end{document}